\documentclass[12pt]{article}
\usepackage[english]{babel}
\usepackage{amsmath}
\usepackage{graphicx}
\usepackage{etoolbox}
\usepackage{changepage}
\usepackage{titlesec}
\usepackage{xurl}
\usepackage[colorlinks=true, allcolors=blue]{hyperref}
\usepackage[parfill]{parskip}
\usepackage[margin=1in]{geometry}
\usepackage{newtx}
\usepackage{float}
\usepackage[numbers,super]{natbib}

\usepackage{caption}
\titleformat*{\section}{\large\bfseries\filcenter} 

\titleformat*{\subsection}{\normalsize\bfseries} 

\usepackage{booktabs}
\usepackage{tabularx}
\usepackage{enumitem}
\usepackage{array}
\usepackage{pifont} 

\newcommand{\yes}{\ding{51}} 
\newcommand{\no}{\ding{55}}  

\newcommand{\systemblock}[5]{%
  \par\medskip
  \noindent\textbf{\normalsize #1}\par
  \smallskip
  \noindent\textit{Method:} #2\par
  \smallskip
  \noindent\textit{Example implementations:} #3\par
  \medskip
  \begin{tabularx}{\linewidth}{@{}X X@{}}
    \toprule
    \textbf{Pros} & \textbf{Cons} \\
    \midrule
    #4 & #5 \\
    \bottomrule
  \end{tabularx}
  \medskip
}

\title{\large \textbf{Tooling for digital accessibility in mathematics:\\ Quickly build compliant course websites that benefit all students}} 
\author{
\normalsize Matthew McMillan, Eli Boyden\\
\normalsize Center for Applied Mathematics, University of Virginia}
\date{} 

\makeatletter 
\patchcmd{\@maketitle}{\begin{center}}{\begin{adjustwidth}{0.5in}{0.5in}\begin{center}}{}{}
\patchcmd{\@maketitle}{\end{center}}{\end{center}\end{adjustwidth}}{}{}
\makeatother

\begin{document}
\raggedright
\maketitle
\thispagestyle{empty}
\pagestyle{empty}


\section*{Abstract}

Public universities in the US are now required to meet digital accessibility (DA) standards under the 2024 updates to Title II of the Americans with Disabilities Act (ADA). For mathematics instructors, this means course materials must be parsable by screen readers, but conventional LaTeX-to-PDF workflows historically do not provide such materials. Despite the availability of Mathematical Markup Language (MathML) as a web standard for accessible math content, supported by all modern browsers, instructor adoption of DA-compliant materials remains very low. This creates a gap between available technology and classroom practice regarding compliance with digital accessibility standards.

This paper makes three contributions toward closing this gap. First, we present a taxonomy of existing approaches for DA-compliant math content, organized by whether they are optimized for print (PDF) or web (HTML) outputs, and we analyze their respective tradeoffs for instructor adoption. Second, we describe and document a proposed implementation: a free software workflow using Obsidian (Markdown-based document authoring and content management tool), Quartz (static site generator), Git (collaboration and version control), and Cloudflare Pages (build and host static sites). This workflow enables math instructors to create, manage, and publish DA-compliant course websites offering MathML from documents that encode math expressions in a TeX-based syntax. The system requires a one-time setup of approximately 1-2 hours, and thereafter site updates occur in minutes, in the cloud, after content authors run a single command. A setup tutorial has been made publicly available by the authors. Third, we present an empirical study of student outcomes across 31 sections of Calculus~II over 6 semesters. Student performance data shows that sections using the proposed course website system outperformed control sections, with the treatment group reaching 2.46 standard deviations above the control mean in the final semester of the study. Although all treatment sections were taught by a single instructor, several lines of evidence, such as the acclimation trajectories of other new instructors on the same course, indicate that the system itself may be a meaningful contributor to the observed performance gains. A survey of student experience shows no statistically significant difference between groups, suggesting that the system does not negatively affect student experience. A proposed second phase of the study will assess barriers to instructor adoption at other institutions.

\section*{Introduction}

\subsection*{Motivation}

On April 24th, 2024 the Department of Justice (DOJ) released its final updates for Title II of the Americans with Disabilities Act (ADA), focusing on ensuring accessible web content in state and local government bodies including public universities and private universities that receive federal funding.\cite{FactSheetNew2024} The technical requirements for this act are dictated by Web Content Accessibility Guidelines (WCAG) 2.1 AA standards produced by the World Wide Web Consortium (W3C). For large public universities, the deadline for compliance is April 2027.

Math instructors are primarily affected by WCAG standard 1.1.1 which requires text alternatives for non-text content.\cite{WebContentAccessibility} W3C recommends Mathematical Markup Language (MathML) for accessible math content, citing its ability to be parsed into audio, braille, and other output formats.\cite{MathematicalMarkupLanguage} The National Center on Accessible Digital Educational Materials and Instruction likewise recommends MathML for its compatibility with screen readers.\cite{BrailleAccessibilityState} Modern screen readers such as Fusion and Job Access With Speech (JAWS) support MathML parsing with audio and braille output.\cite{AccessingMathContent} Empirical evidence also suggests that MathML can improve student and teacher experience compared with traditional textbooks.\cite{lewisUsingAccessibleMath2010}

Despite these standards, recommendations, and technologies, adoption of accessible math content remains low. A survey of 139 physics websites at US post-secondary institutions found accessibility errors on every site.\cite{scanlonPhysicsWebpagesCreate2021} Similar conclusions have been reached regarding digital accessibility in computer science education.\cite{mckieMappingLandscapeDigital2024} MathML adoption specifically has been hampered by historically inconsistent browser support, leading many web developers to offer math notation as rasterized images rather than structured markup.\cite{whiteAccessibilityMathematicalNotation2020} As of 2026, however, MathML is fully supported in every major browser,\cite{MathMLMDN2025} so the remaining barrier is not technological but practical: instructors lack clear, efficient, and well-documented workflows for producing and serving MathML-based course materials.

This paper addresses that practical barrier. We propose and study a complete system in which instructors author and manage their content pages with familiar TeX-based syntax, publish to a website serving MathML, and rebuild the site after content updates by running a single command.

\subsection*{Contributions}

This paper makes three primary contributions:

\begin{enumerate}[noitemsep]
    \item \textbf{Taxonomy of approaches to DA compliance.} We present a structured comparison of existing methods for producing accessible math content, organized by whether they target print (Type~P) or web (Type~W) output, with comments on their tradeoffs for instructor adoption. This taxonomy provides a framework for discussing and evaluating current and future solutions.
    \item \textbf{A documented, easy-to-use workflow.} We describe a specific Type~W implementation based on Obsidian and Quartz, and we provide a step-by-step tutorial designed for instructors with no web development background. We report on the results of using this system over three semesters.
    \item \textbf{Empirical study of student outcomes.} We present student performance data and student experience survey results from 31 sections of Calculus~II over 6 semesters, comparing sections that used the proposed system with those that did not. Although the experimental design has important limitations (most notably that all treatment sections were taught by a single instructor), within-study comparisons help to isolate the contribution of the system from confounding factors.
\end{enumerate}

A proposed second phase, summarized at the end of this article, will study barriers to adoption of this solution by instructors at other institutions.

\subsection*{Design requirements}

To frame the comparison of existing approaches and motivate our design choices, we identify the following theoretical desiderata for DA compliance solutions that could realistically achieve broad adoption among math instructors. A prerequisite for any solution considered here is MathML output: the system must be able to produce MathML embedded in HTML, the recognized standard for accessible math on the web. Since every solution type discussed in this article can produce MathML, we treat it as a baseline capability rather than a distinguishing desideratum, and we compare solutions against the following criteria:

\begin{enumerate}[noitemsep]
    \item \textit{Familiar syntax.} Authors can write math using TeX-based syntax, minimizing retraining.
    \item \textit{Low-friction revision.} Instructors frequently update course materials mid-semester, so the edit-to-published cycle must be simple and fast.
    \item \textit{Minimal technical prerequisites.} The system should not require command-line interactions after the initial setup, nor any significant web-development expertise.
    \item \textit{Free or very low cost.} Budget constraints should not be a barrier to adoption.
    \item \textit{Access management.} Authors should have control over (i) what content is published, and (ii) who can access the site.
    \item \textit{Content management and possession.} The system should support ``Content Management System'' (CMS) features useful for courseware: question banks, transclusion (embed-ability), batch publishing (e.g.\ for homework solutions), links between pages. Authors should also retain possession of their content as portable source files.
    \item \textit{AI workflow friendly.} Content source should reside in the instructor's possession in plain-text formats (e.g.\ Markdown or LaTeX) that integrate readily with AI-assisted workflows for drafting, revising, converting, and batch-editing materials.
\end{enumerate}

These requirements guide the taxonomy discussion and the design of our setup. Note that some of these features go beyond the basic requirements for a successful DA-compliance solution, but they make a system more appealing to instructors and they distinguish our system from certain alternatives.

\section*{Related Work and Taxonomy of Solutions}

\subsection*{Prior scholarship}

Several research efforts have addressed the problem of accessible math content in scholarly articles. Custom systems proposed in the literature include a multimedia platform using a C++ compiler to parse a specialized TeX-like syntax and generate PHP/JavaScript pages,\cite{mackowskiMultimediaPlatformMathematics2018} but resources are not publicly available for other instructors to set up and use this system in their own classrooms. Another system uses Python's SymPy library to convert math input to MathML, parses the MathML into semantic components, and offers interactions allowing users to navigate the mathematical structure of the expressions through a keyboard and speech feedback.\cite{aliAccessibleInteractiveLearning2024} Other notable software tools addressed in research articles include MathPlayer, which provides MathML browser support now common in modern browsers;\cite{soifferMathPlayerV21Webbased2007} MATHTYPE, an editor for MathML and LaTeX content that integrates with other software tools like Word or Canvas;\cite{hindinMATHTYPE662010} and EuroMath, a web-based platform for mathematical communication that appears to be no longer supported.\cite{fitzpatrickEuroMathWebBasedPlatform2020} Recent scholarship has also explored improving PDF accessibility through Markdown conversion using transformer models\cite{duanLayoutAwareTextEditing2026} and by providing supplementary accessible files alongside PDFs.\cite{mittelbachMathMLOtherXML2025}

\subsection*{Taxonomy of solutions}

Our project involves one of various possible approaches to DA-compliance. To contextualize our design choices, we present a taxonomy of the solution space according to the output format primarily targeted, along with examples of software implementations in these types.

The example lists are far from exhaustive. Moreover, there are other solution types that we do not discuss, although they may serve well in individual cases. For example, Microsoft Word documents allow authors to enter equations parsable by screen readers directly from the .docx file. For another example, MathCad is a proprietary system for documentation of engineering calculations that allows users to generate HTML+MathML pages representing the users' notebooks. (SMath is a comparable system with a free version.) However, most mathematics professors do not use these systems because they are already accustomed to the TeX-based syntax which is ubiquitous in mathematics (see Requirement 1), and which is more efficient for writing in large volume than palette-based editors such as in Word or MathCad.

The typesetting and rendering procedures we discuss may be divided according to whether they target PDF or HTML as the primary output. The PDF specification is designed for print: it provides fixed-layout pages well suited to publication for the purposes of reference and archiving, and the LaTeX-to-PDF workflow is deeply established among mathematicians. However, educators have different needs from researchers. They update and rewrite materials frequently, they benefit from modular, internally connected, and easily navigable content, and they serve students who read on devices of various sizes. HTML is a digital-first specification designed for content that is interconnected, reflowable, and dynamically maintained, making it intrinsically better suited for course materials. From the standpoint of digital accessibility, PDF files containing math historically have not been usable by screen readers, whereas HTML with embedded MathML is the recognized accessible format.

The landscape is rapidly evolving. Several significant innovation milestones have been reached in the area of accessible PDFs since this project began and this article was drafted and reviewed. In March of 2026, the international PDF Association announced new standards for accessible math in PDF files in the new PDF/UA-2.\cite{AccessibleMathPDF} Efforts are ongoing to implement updates to LaTeX to allow for generation of PDFs containing accessible math.\cite{LaTeXTaggedPDF} At the same time, there is a growing capacity to make use of these PDFs: although most browsers and PDF readers cannot yet handle embedded math, NVDA and JAWS are now able to do so in Firefox. Converting an old LaTeX file is not necessarily trivial, though, since very many LaTeX packages are partially or completely incompatible with current methods of compiling to the new PDF/UA-2. (See the package compatibility status page.\cite{LaTeXPackageClass}) In any case, instructors desiring a PDF-first approach for their course materials, after considering the other advantages of targeting MathML as described in this paper, are probably well-served to focus on that conversion and the requisite adoption of LuaLaTeX. It is true, of course, that many LaTeX packages are also incompatible with MathJax and KaTeX; our proposal does not aim to maximize compatibility with existing materials.

We denote by ``Type~P'' those systems designed for printable PDF output, and by ``Type~W'' those designed for web HTML output. For DA-compliance, Type~P systems preserve existing LaTeX workflows but require post-conversion to HTML, which introduces friction into the revision cycle. Type~W systems target HTML directly, offering better revision efficiency and higher quality web output, but requiring authors to adopt new (if simpler) writing tools.

It is important to note that DA regulations require the primary learning resource to be accessible. Offering a lower-quality version of content in an accessible form to a student with accessibility needs, while the other students are offered a higher-quality (primary) version in an inaccessible medium, would not be compliant with DA regulations. A system designed to produce PDFs may suffer from ``PDF primacy temptation,'' which is the result of authors focusing on the PDF output quality for most students, while neglecting the quality of an ``accessory'' output that is produced only for compliance purposes. Systems designed for PDF output will tend to contribute DA compliance risk in this way.

The following taxonomy further divides Type P systems into P1 and P2, and Type W systems into W1, W2, and W3:

\setlist[itemize]{noitemsep, topsep=2pt, leftmargin=*}

\systemblock
  {P1 --- Machine conversion: PDF $\rightarrow$ MathML}
  {Software tools, typically implementing ML character recognition, possibly with LLM support, convert the PDF file itself to HTML webpages with embedded MathML. These files are then loaded in the public pages of a course website. Attention is required for links between files.}
  {MathPix (math-aware OCR); EditTrans (LLM-based), PDFix}
  {
    \begin{itemize}
      \item Meets desiderata 1, 3, 4, 7
      \item Preserve document shape and layout
      \item Continue using LaTeX source
      \item Same method applies to old materials
    \end{itemize}
  }
  {
    \begin{itemize}
      \item Fails desiderata 2, 5, 6
      \item Poor author control over resulting HTML
      \item Output HTML requires manual cleanup
      \item Inefficient cycle for continuous revision
      \item DA risk: PDF primacy temptation
    \end{itemize}
  }

\systemblock
  {P2 --- Machine conversion: LaTeX $\rightarrow$ MathML}
  {Software tools convert LaTeX source code to HTML with MathML.}
  {LaTeXML (used by arXiv.org); Pandoc; tex4ht; Math-Core}
  {
    \begin{itemize}
      \item Meets desiderata 1, 3, 4, 7
      \item Preserve document structure
      \item Continue using LaTeX source
      \item Same method applies to old materials
    \end{itemize}
  }
  {
    \begin{itemize}
      \item Fails desiderata 2, 5, 6
      \item Moderate control over resulting HTML
      \item Must adjust syntax, settings each file
      \item Conversion time discourages frequent revision
      \item DA risk: PDF primacy temptation
    \end{itemize}
  }


\systemblock
  {W1 --- Create and manage content in an LMS}
  {Learning Management Systems like Canvas or Blackboard allow instructors to add math expressions to internal content pages using a built-in equation editor, producing output compatible with built-in accessibility software functions.}
  {Canvas, Blackboard}
  {
    \begin{itemize}
      \item Meets desiderata 2, 3, 4, 5
      \item Some limited CMS functionality
      \item Good control over HTML (less over MathML)
      \item Institutional site access control is built in
    \end{itemize}
  }
  {
    \begin{itemize}
      \item Fails desiderata 1, 6, 7
      \item Layout and TeX limitations (vs.\ LaTeX)
      \item Poor math editor UI, inadequate for heavy use
      \item Cannot export math source
      \item Cannot export PDF versions
      \item Doesn't apply to old materials
    \end{itemize}
  }

\systemblock
  {W2 --- Create and manage Markdown+TeX source in pro web CMS}
  {Web publishing systems such as Wordpress or Drupal can serve math using HTML including MathML. However, these systems are designed primarily for general web publishing rather than courseware authoring. They require substantial web development skills and resources to maintain, and they carry significant subscription or licensing costs.}
  {Wordpress, Drupal, SquareSpace, WIX}
  {
    \begin{itemize}
      \item Meets desiderata 1, 5, 6
      \item Myriad website features are possible
    \end{itemize}
  }
  {
    \begin{itemize}
      \item Fails desiderata 2, 3, 4, 7
      \item Some layout and TeX limitations (vs.\ LaTeX)
      \item Requires web dev skills: setup and ongoing
      \item CMS features for courseware aren't bundled
    \end{itemize}
  }

\systemblock
  {W3 --- Create and manage Markdown+TeX source in personal local CMS}
  {Write course material in local files using Markdown and a TeX variant. Use a static site generator tool to produce HTML with MathML.}
  {Obsidian, LogSeq, Jupyter Notebooks, Hugo, Notion, Craft; Jamstack generators}
  {
    \begin{itemize}
      \item Meets desiderata 1, 2, 3, 4, 5, 6, 7
      \item Easy set up compared with pro CMS
      \item Expedite authoring with plugins
      \item CMS functionality fits courseware needs
    \end{itemize}
  }
  {
    \begin{itemize}
      \item Fails desiderata [none]
      \item Some layout and TeX limitations (vs.\ LaTeX)
      \item PDF output is secondary, lower quality
    \end{itemize}
  }

\par\bigskip
\noindent\textbf{Summary}\par
\medskip

\begin{tabularx}{\linewidth}{@{}l X X@{}}
\toprule
\textbf{System} & \textbf{Primary Strength} & \textbf{Primary Weaknesses} \\
\midrule
P1 & Simplest DA solution & Low quality HTML output; opaque and slow conversion discourages revision \\
P2 & Robust LaTeX-first & Setup time per file; opaque and slow conversion; ``PDF primacy'' temptation \\
W1 & Easy start & Slow to use; content trapped; no PDFs \\
W2 & Myriad website features possible & Prohibitive prereq skills \\
W3 & Web-first yet serves PDF; easy to write; compelling CMS features & Mildly technical one-time setup; reduced PDF layout control vs.\ LaTeX \\
\bottomrule
\end{tabularx}

\setlist[itemize]{}

\medskip

Evaluated against the design requirements identified in the introduction, Type~W3 systems best balance the competing demands of accessibility compliance, low to moderate startup overhead, content management functionality that is useful for courseware, and revision cycle efficiency. Table~\ref{tab:scorecard} summarizes this evaluation as a scorecard of the five system types against the numbered design requirements. Our implementation provided below falls into the W3 category.

\begin{table}[htb!]
    \centering
    \begin{tabularx}{\linewidth}{@{}X c c c c c@{}}
        \toprule
        \textbf{Design requirement}
        & \shortstack{\textbf{P1}\\[1pt] \scriptsize PDF$\rightarrow$ML}
        & \shortstack{\textbf{P2}\\[1pt] \scriptsize LaTeX$\rightarrow$ML}
        & \shortstack{\textbf{W1}\\[1pt] \scriptsize LMS}
        & \shortstack{\textbf{W2}\\[1pt] \scriptsize Pro CMS}
        & \shortstack{\textbf{W3}\\[1pt] \scriptsize Ours} \\
        \midrule
        1.~Familiar TeX syntax                & \yes & \yes & \no  & \yes & \yes \\
        2.~Low-friction revision              & \no  & \no  & \yes & \no  & \yes \\
        3.~Low expertise                      & \yes & \yes & \yes & \no  & \yes \\
        4.~Low cost                           & \yes & \yes & \yes & \no  & \yes \\
        5.~Access management                  & \no  & \no  & \yes & \yes & \yes \\
        6.~Content management and possession  & \no  & \no  & \no  & \yes & \yes \\
        7.~AI workflows friendly              & \yes & \yes & \no  & \no  & \yes \\
        \bottomrule
    \end{tabularx}
    \caption{Scorecard of MathML solution types against the design requirements. The numbered rows correspond to the distinguishing desiderata listed in the introduction.}
    \label{tab:scorecard}
\end{table}

\section*{System Overview}

This section describes the system we have implemented. We first present the daily instructor workflow to convey the simplicity of ongoing use, then describe the one-time setup process. For a complete step-by-step tutorial designed for instructors with no web development background, see our public tutorial at \texttt{\href{https://quartz-mathml-tutorial.pages.dev/}{quartz-mathml-tutorial.pages.dev}}. Initial setup is expected to take 1-2 hours following the tutorial; after setup, authors work entirely within Obsidian to write content and initiate site rebuilds. The site rebuild is triggered by a command in Obsidian and takes place in the cloud. It requires 2-5 minutes to build our site which hosts all material for six courses.

\subsection*{Daily workflow}
Content is written in Obsidian, a highly extensible note-taking application that views and manages a folder, referred to as a ``vault,'' containing Markdown files stored locally on the instructor's computer. Math is written in a TeX-like syntax (with single \$ for inline math, double \$\$ for display math). Obsidian provides a streamlined writing interface with features such as LaTeX autocompletion, internal linking between pages, and transclusion (embedding a part of one page's contents inside another page). These features enable modular courseware: instructors can maintain a personal library of theory notes, examples, and exercises, transclude them into semester-specific pages, manage question banks with linked solutions, and batch-publish homework solutions by toggling a ``property'' name in the file frontmatter.

Publishing follows the ``commit and sync'' command provided by the Obsidian Git plugin. This command automatically synchronizes the local folder with the cloud GitHub repository using a ``commit'' to save and label the current state. Cloudflare Pages detects the update and automatically rebuilds the website in the cloud. The author does not need to take any further steps, and the new site will be live in a few minutes.

\subsection*{One-time setup}

\begin{figure}[H]
    \centering
    \includegraphics[width=0.55\linewidth]{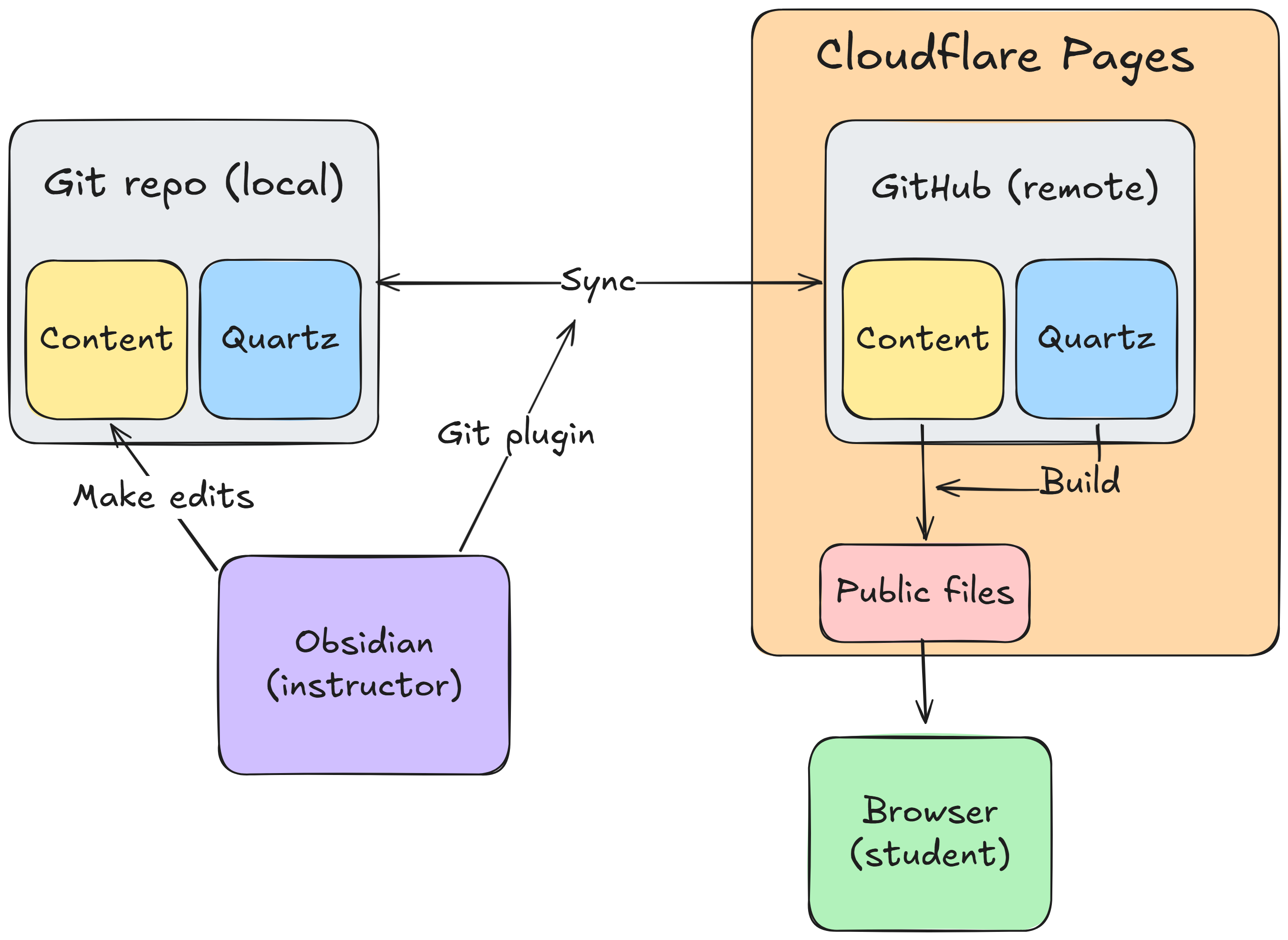}
    \caption{System workflow. The instructor interacts only with Obsidian (left). A single ``Commit and Sync'' command triggers automatic building and deployment of the course website (right).}
    \label{fig:sys-proc}
\end{figure}

The system uses four applications (see Figure~\ref{fig:sys-proc}):

\begin{enumerate}[noitemsep]
    \item \textbf{Obsidian} --- Composition and content management. Obsidian's simple Markdown file format, internal link management, and vast free plugin ecosystem (including \textit{Obsidian Git} for version control) make it well-adapted for courseware managed by individual math instructors. Although a setup without plugins (besides \textit{Obsidian Git}) would suffice, we recommend considering use of the following plugins while optimizing a workflow over time: \textit{Admonition}, \textit{Advanced Cursors}, \textit{Advanced Tables}, \textit{Better Export PDF}, \textit{Copy Block Link}, \textit{Easy Bake}, \textit{Excalidraw}, \textit{File Cooker}, \textit{Git}, \textit{Linter}, \textit{Multi Properties}, \textit{Note Refactor}, \textit{PDF Break Page}, \textit{Quick Latex} (or \textit{Latex Suite}), \textit{Quick Switcher++}, \textit{Recent Files}, \textit{Settings Search}, \textit{Style Settings}, \textit{Table of Contents}, \textit{Tag Wrangler}, \textit{Text Format}, \textit{Text Transporter}, \textit{TikZJax}, \textit{Underline}.
    \item \textbf{Git/GitHub} --- Version control and cloud storage. Content files are synced with a GitHub repository, enabling multi-instructor collaboration and saved snapshots, and providing the connection point for automated deployment of the site.
    \item \textbf{Quartz} --- Static site generator tool. Quartz is designed and maintained by Jacky Zhao (\texttt{\href{https://jzhao.xyz/}{jzhao.xyz}}) specifically to build websites from Markdown collections in Obsidian or similar apps. It transforms Markdown files into a collection of HTML, CSS, and JavaScript files. Critically, Quartz supports MathML generation with minor customization that we explain in the tutorial. The simplest method just passes an option to the KaTeX parsing engine. We also developed a custom Quartz plugin using the Temml library (via the open-source rehype-mathml project); this provides wider TeX syntax coverage and some MathML rendering improvements. Note that Quartz requires NodeJS as a dependency.
    \item \textbf{Cloudflare Pages} --- Free web hosting with automatic build and deployment triggered by changes to the GitHub repository. Other hosting options (Vercel, Netlify, self-hosting) are also possible.
\end{enumerate}

Customization is possible at each layer. For example, Quartz supports custom SCSS for styling; we use this to create custom ``callout'' block styling for concept definitions, examples, exercises, and solutions as are commonly included in math course materials. These customizations and all other setup details are covered in the tutorial.

\subsection*{Documentation links and website screenshots}

\begin{figure}[htb!]
    \centering
    \includegraphics[width=0.9\linewidth]{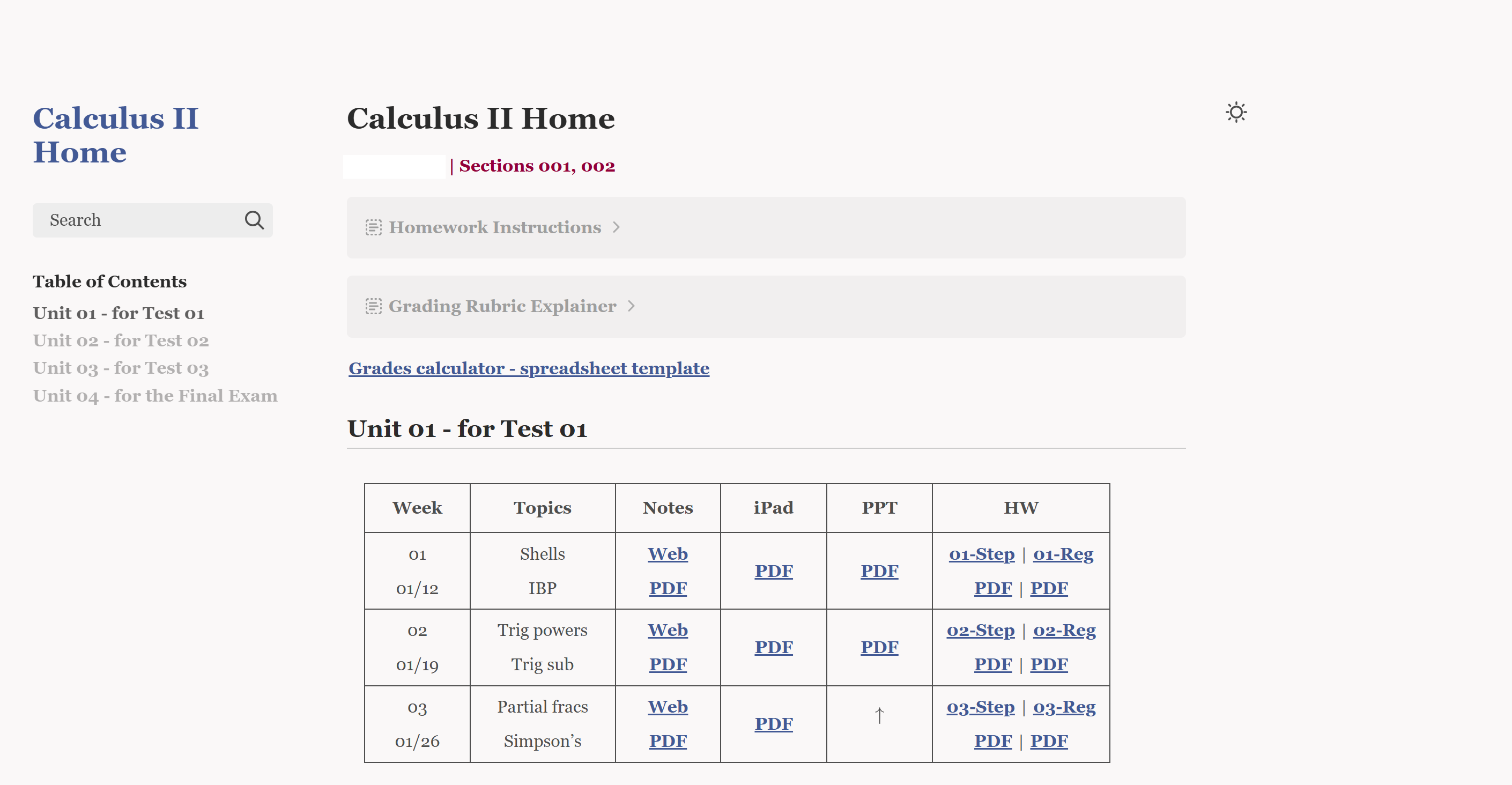}
    \caption{Course website: Calculus II home page with links to materials pages arranged by week. Tables are easy to manage in Obsidian.}
    \label{fig:WebsiteHome01}
\end{figure}

\begin{figure}[htb!]
    \centering
    \includegraphics[width=0.9\linewidth]{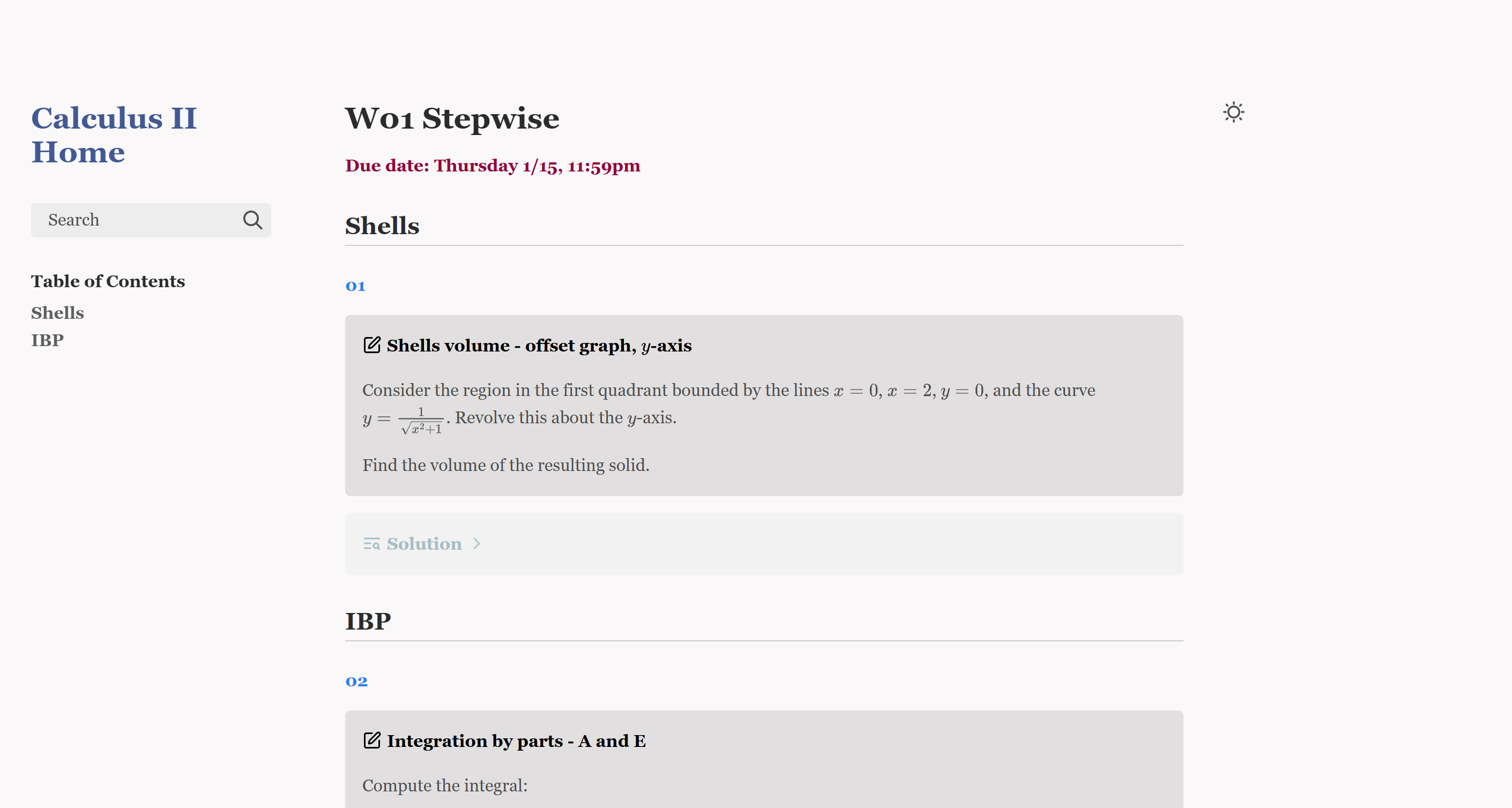}
    \caption{Course website: Sample homework page showing expandable solution callouts. Solutions displayed in the callouts on a given homework page are easily published or hidden in a batch.}
    \label{fig:WebsiteHW01}
\end{figure}

See Figure~\ref{fig:WebsiteHome01} for a sample of the course homepage, and Figure~\ref{fig:WebsiteHW01} for a sample of a homework assignment. Documentation for each software tool is currently found at the following pages:

Obsidian: \href{https://help.obsidian.md/}{https://help.obsidian.md/}\\
Git: \href{https://git-scm.com/docs}{https://git-scm.com/docs}\\
GitHub: \href{https://docs.github.com/en}{https://docs.github.com/en}\\
Quartz: \href{https://quartz.jzhao.xyz/}{https://quartz.jzhao.xyz/}\\
NodeJS: \href{https://nodejs.org/docs/latest/api/}{https://nodejs.org/docs/latest/api/}\\
Cloudflare: \href{https://developers.cloudflare.com/directory/}{https://developers.cloudflare.com/directory/}

\subsection*{Variants and alternatives}

Variants on our solution within Type~W3 are also reasonable. For example, instead of using Cloudflare Pages, one can use GitHub Pages for site hosting as well as for the content repository. We chose Cloudflare Pages because it does not expose all source files to the public. At present a GitHub repository is either completely public or private, and only public repositories are publishable through GitHub Pages. By keeping our GitHub repository private, we can include in the one main repository any homework solutions, exams, and other contents we do not wish to make public at all times. It is also possible to set up password protection for a site for free using a Cloudflare Worker, although the setup for this is more complex.

For another variant system, one could use Obsidian only to create the content, and then export each file to HTML (via an Obsidian plugin) and load the HTML files directly into an LMS such as Canvas. This would provide institutional access control through the LMS, but it would forego much of the content management functionality and the efficiency of automated publishing.

An important motivation for our system, and one reason we believe it could improve adoption of DA-compliant materials, is that it holds promise to benefit \textit{all} students --- not just those who rely on accessibility features. Ease of accessing and navigating a well-organized course website may translate to more frequent and efficient study for all students, and ease of managing content may reduce the friction experienced by instructors making continual improvements to materials. We examine evidence for these hypotheses in the next section.

\section*{Student Impact: Evidence from Calculus II}

This section presents an empirical study assessing the impact on students of the system described above. We examine student impact in two dimensions: student performance (assessed through overall course grades) and student experience (assessed through a survey). We emphasize that this study did not attempt to assess the \textit{degree} of DA compliance of the generated course materials, nor the experience of low-vision students specifically. (We consider MathML as a proxy for DA compliance of math course materials; other aspects of DA standards are not specific to math.) The goal of this empirical study was to determine whether use of the system was associated with a positive, negative, or indeterminate impact on the general student population.

\subsection*{Setup}

Our study involved a multi-section Calculus~II course over 6 semesters, ``S23'' (Spring 2023) through ``F25'' (Fall 2025), in 31 distinct sections. An average of 197 students per semester (range: 145-272) enrolled in 4-7 sections, with section sizes averaging 38 students (range: 18-48, smallest three being 18, 25, 29). All 6 instructors were full-time permanent teaching staff, either senior lecturers or teaching-track professors (none were grad students, postdocs, or adjunct).

All students in all sections took the same four exams per semester (three midterms and one cumulative final). These exams determined 70-75\% of each student's grade. Exams were conducted in person without technology access. Exams were drafted by one instructor in S23-S24 and by a different instructor in F24-F25; both of these instructors had taught this course many times. Homework problems were mostly the same across sections, and included both Cengage WebAssign problems and handwritten work. Instructors differed in lecture style (at least one used a flipped classroom), lecture notes, use of class time, and experience with the course.

We divide the 31 sections into ``control'' and ``treatment'' groups. Control sections (26 total) did not in any way use the course website we described, instead they offered PDF files through Canvas for lecture slides, in-class worksheets, and solutions. Treatment sections (5 total, starting F24) were taught by a single instructor who used the system described in this article in a phased rollout that we explain now. This phased rollout is important for interpreting our data.

In F24, the treatment instructor used Obsidian to type Markdown+TeX files for course materials for two sections (i.e.\ for lecture notes and written homework problems), and then PDFs were generated from these source files and offered to students inside Canvas, the institutional LMS (the same delivery format used by control sections). The F24 semester therefore serves as a same-instructor baseline for the subsequent treatment semesters: the instructor used the proposed authoring tools but did not yet provide students with the content organized on webpages. Note, however, that F24 was also the treatment instructor's first semester on this particular course.

In S25, there was one treatment section (with 31-34 students) and those students received all materials on a course website produced from the Markdown+TeX source files, although homework solutions were still provided as scanned PDFs. In F25, the (two) treatment sections received all materials in a website generated by Quartz (including solutions available directly below problems, after submission deadlines, as shown in Figure~\ref{fig:WebsiteHW01}) using the full system as described in the System Overview above.

One additional institutional change must be noted. Between S25 and F25, a policy change required students to hold formal college credit for Calculus~I before enrolling in Calculus~II. This change might be expected to improve the preparation of students equally in both treatment and control groups in F25. We return to this point in the discussion.

Data on student performance was obtained from the internal Office of Institutional Research \& Analytics (IRA). See Table~\ref{tab:grade_distributions_full} in the Appendix for complete source data as used in our analysis. For each section, we gathered and computed: (i)~the percentage of students receiving A-band grades (A+, A, A$-$), (ii)~the percentage receiving DFW, and (iii)~the GPA as an average of standard grade-point values, averaged over students in the section.

Data on student experience was collected through a survey administered by a central university office to the F25 cohort at the end of the semester. The survey drew 55 responses from treatment students (out of 89 in that group) and 16 from control students (also out of 87 in that group). The central office administers such surveys as part of a broad study authorized by the IRB at this university.

\subsection*{Performance results}

\begin{figure}[htb!]
    \centering
    \includegraphics[width=1\linewidth]{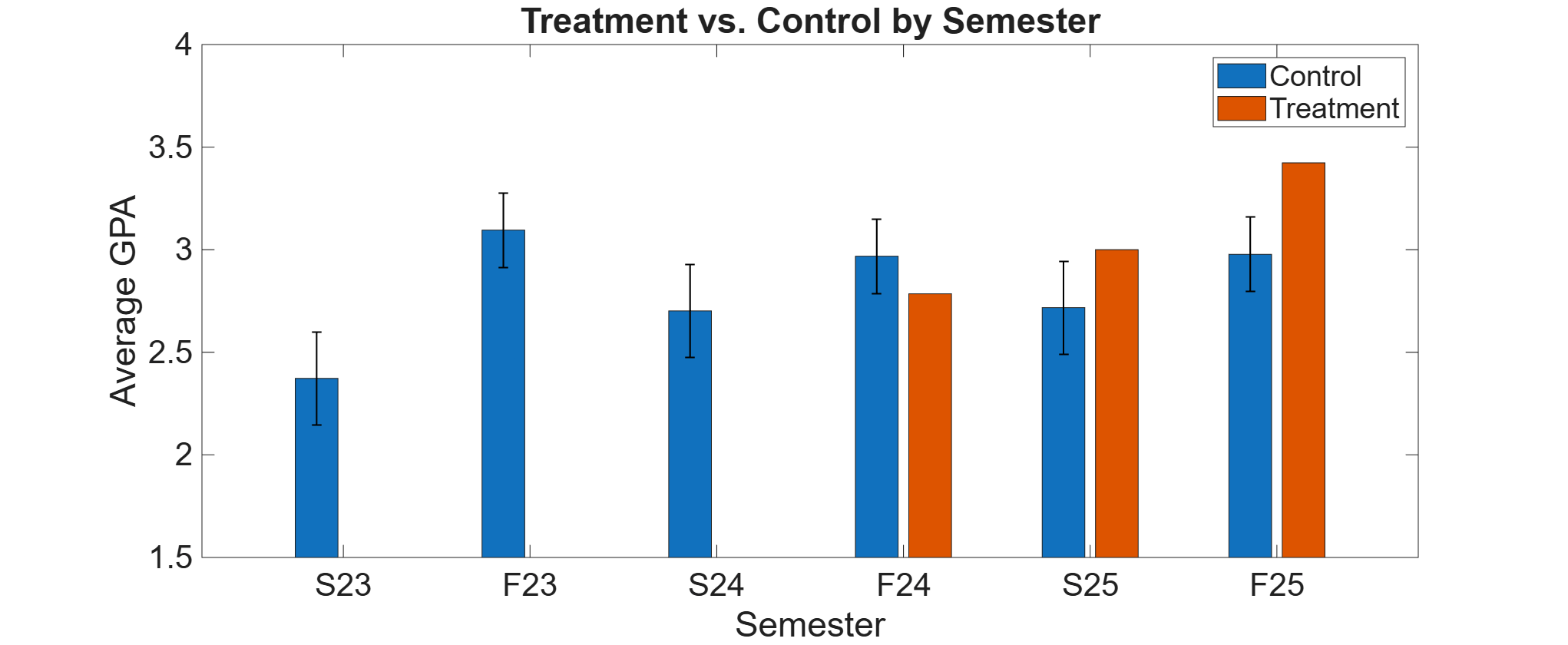}
    \caption{GPA by semester. Blue bars: control sections; red bars: treatment sections. Error bars represent $\sigma_\text{fall}$ or $\sigma_\text{spring}$, the population standard deviation of control section GPAs by season. These error bars indicate the typical spread of section GPAs around the seasonal control average.}
    \label{fig:GPA_semesters}
\end{figure}

In Figure~\ref{fig:GPA_semesters} we depict the section GPA averaged across sections in each semester. Blue bars represent control section averages; red bars represent treatment section averages. By semester, the control bars account for 4, 6, 5, 5, 4, and 2 sections; the treatment bars account for 2, 1, and 2 sections (starting F24).

The error bars on control sections represent the population standard deviation of all control sections in the corresponding season ($\sigma_\text{fall}$ computed over 13 fall sections, $\sigma_\text{spring}$ over 13 spring sections). We separate seasons because we observe a stable pattern of stronger performance in fall cohorts than in spring cohorts, which is visible in the oscillating blue bars.

Three features of the data are notable:

\begin{enumerate}[noitemsep]
    \item The F24 treatment bar (same-instructor baseline with PDF delivery) sits approximately $\sigma_\text{fall}$ \textit{below} the control average, consistent with the instructor's first semester on the course.
    \item The S25 treatment bar (first semester using a course website) exceeds the control average by more than $\sigma_\text{spring}$, and is higher than the F24 treatment bar despite the typical pattern of \textit{lower} grades in spring semesters.
    \item The F25 treatment bar (full system operational) is 2.46 increments of $\sigma_\text{fall}$ above the F25 control bar, a significant separation.
\end{enumerate}

The jump in performance between F24 and S25 is particularly informative. Since the same instructor taught both semesters using Markdown+TeX source materials, but switched from PDF delivery (F24) to website delivery (S25), the performance improvement between these semesters, and especially the reversal of the usual fall-to-spring decline, is consistent with the course website making a meaningful difference.

\begin{figure}[h]
    \centering
    \includegraphics[width=1\linewidth]{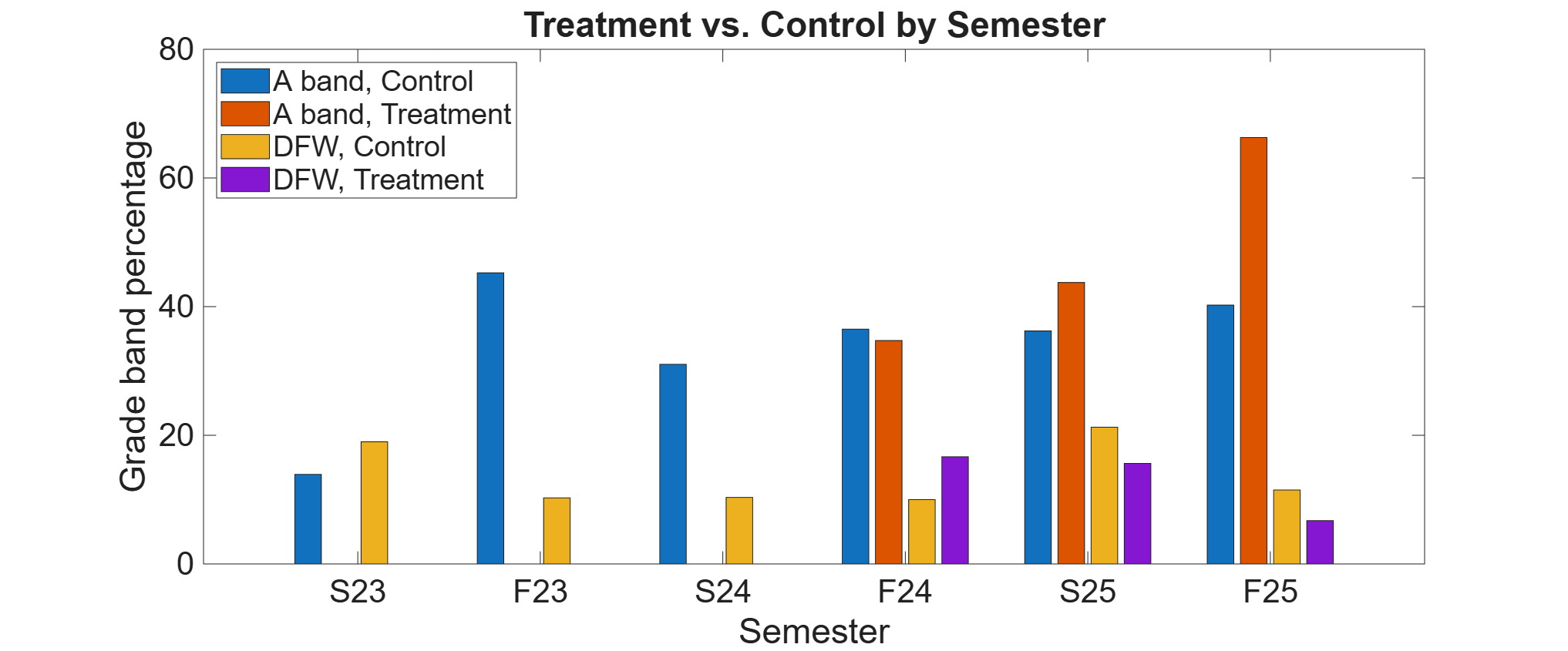}
    \caption{A-band and DFW percentages by semester.}
    \label{fig:ADFW_semesters}
\end{figure}

Figure~\ref{fig:ADFW_semesters} depicts the percentages of students receiving A-band grades and DFW outcomes. The A-band growth in treatment sections is even more pronounced than the GPA trend. The treatment DFW rate also steadily decreases and is lower than the control in S25 and F25, when the course website was in use.

\subsection*{Survey results}

The centrally administered survey included the following questions relevant to this study:

\begin{enumerate}[noitemsep]
\small
    \item[$Q_1:$] ``The primary course materials of my calculus class (notes, slides, videos, worksheets, problems, textbook) offered in my course/section were very helpful overall.''
    \item[$Q_2:$] ``It was easy to understand the course materials when I viewed and studied them the first time through.''
    \item[$Q_3:$] ``It was easy to navigate the course materials while I was studying later on for exams or solving problems.''
    \item[$Q_4:$] ``Ease of navigation of course materials was important for supporting my problem-solving work.''
    \item[$Q_5:$] ``Ease of navigation of course materials was important for supporting my review study for exams.''
\end{enumerate}

Respondents selected one option from the list: ``Strongly agree,'' ``Agree,'' ``Neutral,'' ``Disagree,'' ``Strongly disagree.'' We analyzed responses using two coding schemes. In \textit{collapsed coding}, responses were mapped to three categories reflecting direction of agreement: ``Strongly agree'' and ``Agree'' $\to$ 3, ``Neutral'' $\to$ 2, ``Disagree'' and ``Strongly disagree'' $\to$ 1. This coding focuses on the direction of agreement rather than its intensity, which may be more robust given the small control group size. In \textit{full Likert coding}, responses were mapped to 5 (``Strongly agree'') through 1 (``Strongly disagree'') in order. Results for both codings are reported in Table~\ref{tab:survey_both}.

\begin{table}[htb!]
    \centering
    \begin{tabular}{l c c c c c c}
        \toprule
         & \multicolumn{3}{c}{\textbf{Collapsed (1--3)}} & \multicolumn{3}{c}{\textbf{Full Likert (1--5)}} \\
        \cmidrule(lr){2-4} \cmidrule(lr){5-7}
        \textbf{Q} & \textbf{Treatment} & \textbf{Control} & $\sigma_i$ & \textbf{Treatment} & \textbf{Control} & $\sigma_i$ \\
        \midrule
        $Q_1$ & 2.75 & 2.75 & $\pm 0.52$ & 3.95 & 4.06 & $\pm 0.75$ \\
        $Q_2$ & 2.44 & 2.19 & $\pm 0.73$ & 3.51 & 3.31 & $\pm 0.85$ \\
        $Q_3$ & 2.84 & 2.69 & $\pm 0.43$ & 4.13 & 4.06 & $\pm 0.72$ \\
        $Q_4$ & 2.76 & 2.69 & $\pm 0.49$ & 3.98 & 4.00 & $\pm 0.74$ \\
        $Q_5$ & 2.78 & 2.81 & $\pm 0.41$ & 4.09 & 4.19 & $\pm 0.72$ \\
        \bottomrule
    \end{tabular}
    \caption{Survey results: mean values and pooled standard deviations under both coding schemes. Higher values indicate stronger agreement.}
    \label{tab:survey_both}
\end{table}

\begin{figure}[htb!]
    \centering
    \includegraphics[width=1\linewidth]{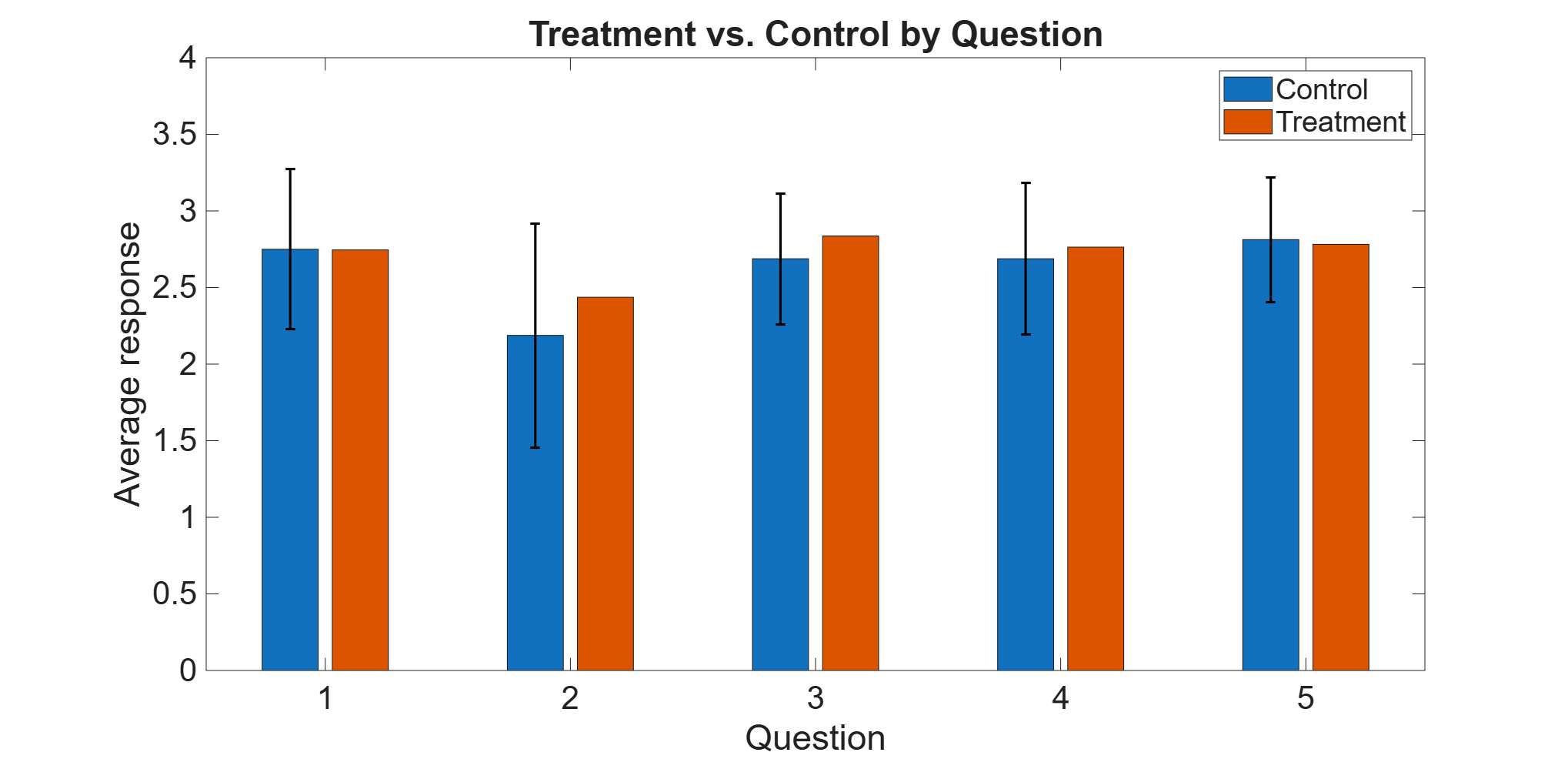}
    \caption{Survey response averages by question (collapsed coding).}
    \label{fig:Survey}
\end{figure}

The results are also depicted in Figure~\ref{fig:Survey} using the collapsed coding. For each question, the error bar represents a pooled standard deviation for that question computed as:
$$
\sigma_{i} \;=\; \sqrt{\frac{1}{55+16}\Big((55)\sigma_{T_i}^2+(16)\sigma_{C_i}^2\Big)}
$$
where $\sigma_{T_i}$ and $\sigma_{C_i}$ are the population standard deviations of the treatment and control groups' coded responses for question $Q_i$.

Under the collapsed coding, the treatment trends slightly higher on $Q_2$ and $Q_3$, but all differences are well within $\frac{1}{2}\,\sigma_i$ of the control. Under the full Likert coding, the treatment lags slightly below the control on $Q_1$ and $Q_4$. A Welch's t-test (one-tailed, unequal variances) yielded $p = 0.15$ for $Q_2$ and $p = 0.18$ for $Q_3$ under the collapsed coding; all other $p$-values were higher. So we do not see a statistically significant difference under either coding scheme.

\subsection*{Limitations}

The experimental design of this study carries several important limitations. We discuss each now alongside available mitigating evidence.

\textit{Single instructor.} All treatment sections were taught by a single instructor who taught no control sections, so instructor effects cannot be fully separated from treatment effects in this design. The within-instructor comparison between F24 (PDF delivery) and S25/F25 (website delivery) provides partial control for the instructor variable, since the same person taught all three semesters using very similar source materials but with different delivery methods. The performance increase between F24 and S25 coincides with the introduction of the course website rather than with a change in treatment instructor.

\begin{figure}[htb!]
    \centering
    \includegraphics[width=1\linewidth]{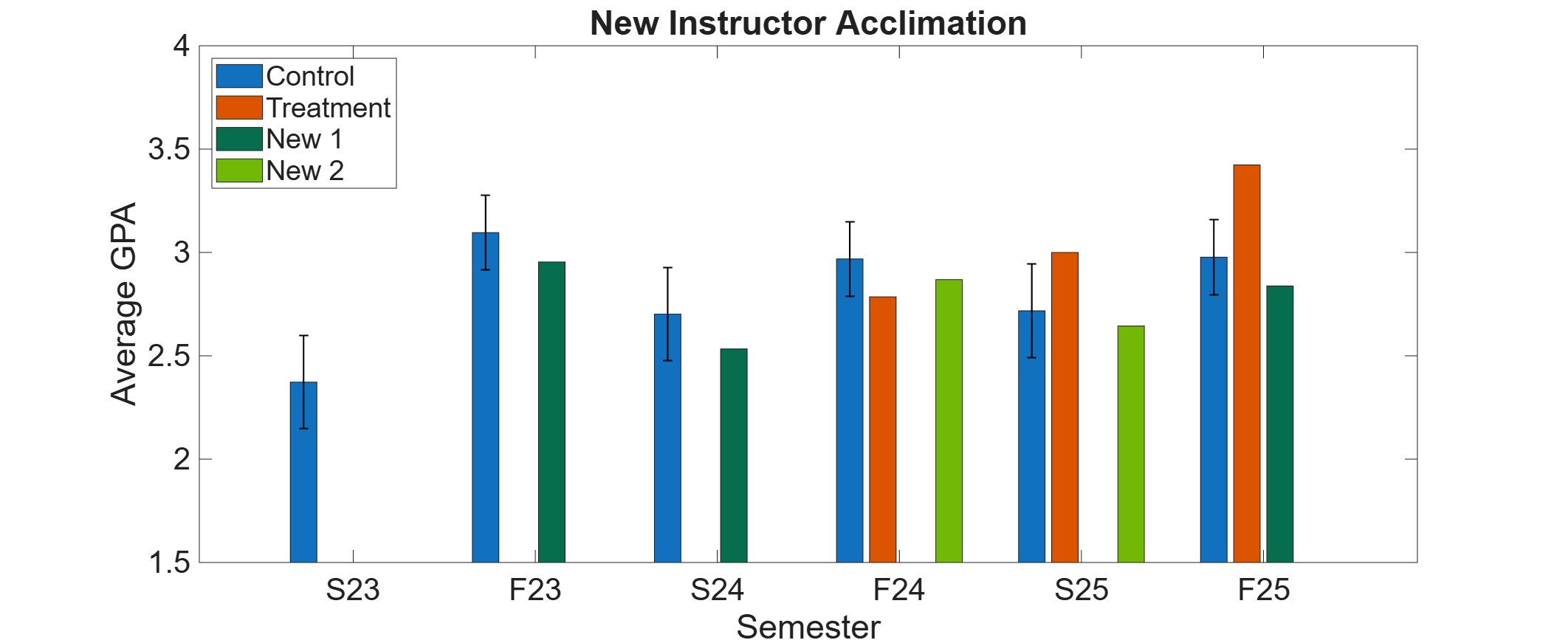}
    \caption{Semester-by-semester performance for two other instructors who were new to the course during the study period and who did not implement the proposed system, alongside Figure~\ref{fig:GPA_semesters} bars. Neither shows gains comparable to the treatment group, relative to the control.}
    \label{fig:new_instructor}
\end{figure}

\textit{Instructor acclimation.} The treatment instructor was new to this course in F24, so some performance improvement in the subsequent semesters may be expected as the instructor acclimated to the course content, learning objectives, and the style and expectations of the common assessments. The course materials were also adjusted in minor ways between F24 and S25 to better align with the common exams, and the use of class time shifted from blackboard lecturing (F24) to annotating projected notes from an iPad (S25, F25). To help assess whether instructor acclimation alone could account for the observed trajectory, Figure~\ref{fig:new_instructor} shows the performance over comparable multi-semester periods for two other instructors who were also new to this course during the study window but who did not implement the system described in this article. (New instructor 1 taught this course only F23, S24, F25, while New instructor 2 taught only F24, S25.) Neither shows performance gains of comparable magnitude relative to the control, suggesting that acclimation to the course is not sufficient on its own to explain the treatment group's trajectory.

\newpage

\textit{Prerequisite policy change.} Between S25 and F25, the institution began requiring formal Calculus~I credit for enrollment in Calculus~II, which may have improved the preparation of F25 students. However, this policy change would be expected to affect both treatment and control groups equally. Inspecting the control bars in Figure~\ref{fig:GPA_semesters}, we observe that the control group GPA does not increase from F24 to F25; indeed, control GPAs are flat or slightly declining over the F23--F24--F25 sequence. This suggests that the prerequisite policy change is not the primary driver of the treatment group's performance gains in F25.

\textit{Small survey control group.} The student experience survey drew only 16 responses from the control group, which significantly limits the statistical power of the survey analysis. Consequently, the null result on the survey should not be interpreted as strong evidence that the system has no effect on student experience, but neither does it indicate a negative effect. Extending the survey to additional cohorts is a goal of future work.

\subsection*{Discussion}

Taken together, the performance data present a consistent picture. The treatment instructor's trajectory shows a marked improvement that begins precisely when the course website is introduced (S25) rather than when the instructor begins teaching the course (F24). This improvement grows further in F25 when the full system is operational including the website access to homework solutions after submission deadlines. The $2.46\,\sigma_\text{fall}$ separation between treatment and control in F25 is notable given that $\sigma_\text{fall}$ reflects the normal range of variation across 6 instructors (13 fall control sections) with varying levels of course experience. A result of this magnitude would be difficult to achieve through instructor individuality or acclimation alone. Two other new instructors on the same course do not show similar gains over comparable time frames, and the prerequisite policy change does not appear to have benefited control sections. While the single-instructor design means we cannot definitively attribute the performance gains to the courseware system, the combination of these lines of evidence supports the conclusion that the system is a meaningful contributor.

The survey results do not show a statistically significant difference between treatment and control groups on any question under either coding scheme. The small control group (16 respondents) limits what may be concluded with any confidence. The data does, however, indicate that the system is unlikely to have a substantial negative effect on student experience. This is a useful finding for the adoption problem, since a negative effect on student experience would weigh against wide adoption of the system.

\section*{Future Work: Study of Adoption Barriers}

Having established a DA-compliant solution and studied its effects on students, we turn to the next major question: what are the barriers to adoption of this solution by other faculty?

We propose the identification of a small group of instructors at other institutions with varied teaching contexts who are willing to adopt the system for their own courses. We will collect data through instructor surveys (before, during, and after adoption) and through records of issues encountered during the setup and maintenance process. The initial survey will assess participants' backgrounds in programming, web development, LaTeX use, and course material creation; the final survey will assess satisfaction after several semesters using our system.

As a secondary objective, participating faculty would contribute non-private student performance data from before and after adopting the system. Aggregating this data across multiple institutions and instructors would strengthen the evidence for the student performance findings reported here by addressing the single-instructor limitation of the current study, and it would provide further motivation for broad adoption.

\section*{Acknowledgments}

We thank Meiqin Li for helpful comments and discussion, and Lindsay Wheeler for helping us source grades data through the UVA IRA. We also thank the UVA Center for Teaching Excellence for supporting work by the second author.

\vspace{4\baselineskip}\vspace{-\parskip} 
\footnotesize 
\bibliographystyle{unsrtnat} 
\bibliography{References}


\section*{Appendix}

\begin{table}[H]
    \centering
    \small
    \begin{tabular}{l c c c c c c c c c c c c}
        \toprule
        \textbf{Term} & \textbf{Sec.} & \textbf{A+} & \textbf{A} & \textbf{A--} & \textbf{B+} & \textbf{B} & \textbf{B--} & \textbf{C+} & \textbf{C} & \textbf{C--} & \textbf{DFW} & \textbf{N} \\
        \midrule
        2023 Spring & 1 &  1 &  2 &  3 &  2 &  0 &  1 &  3 &  2 &  7 & 11 & 32 \\
                    & 2 &  1 &  2 &  0 &  6 &  7 &  3 &  6 &  4 &  6 &  8 & 43 \\
                    & 3 &  1 &  4 &  3 &  7 &  8 &  5 &  1 &  5 &  2 &  7 & 43 \\
                    & 4 &  1 &  1 &  3 &  6 &  3 &  4 &  3 &  5 & 10 &  4 & 40 \\
        \addlinespace
        2023 Fall   & 1 &  3 &  7 & 10 &  3 &  7 &  3 &  2 &  4 &  2 &  3 & 44 \\
                    & 2 &  2 &  8 &  7 &  6 &  6 &  4 &  2 &  3 &  1 &  4 & 43 \\
                    & 3 &  8 &  9 &  7 &  3 &  5 &  4 &  1 &  4 &  0 &  4 & 45 \\
                    & 4 &  6 & 15 &  7 &  5 &  6 &  1 &  0 &  0 &  2 &  4 & 46 \\
                    & 5 &  2 &  9 &  6 &  4 &  7 &  2 &  1 &  1 &  0 &  8 & 40 \\
                    & 6 &  1 &  7 &  5 &  8 & 10 &  3 &  1 &  5 &  1 &  4 & 45 \\
        \addlinespace
        2024 Spring & 1 &  1 &  7 &  4 &  5 &  2 &  3 &  6 &  6 &  0 &  1 & 35 \\
                    & 2 &  1 &  5 &  5 &  2 &  3 &  4 &  4 &  1 &  4 &  3 & 32 \\
                    & 3 &  3 &  2 &  3 &  2 &  3 &  1 &  3 &  4 &  2 &  6 & 29 \\
                    & 4 &  2 &  1 &  4 &  5 &  1 &  3 &  4 &  4 &  5 &  2 & 31 \\
                    & 5 &  1 &  2 &  4 &  0 &  0 &  1 &  1 &  2 &  4 &  3 & 18 \\
        \addlinespace
        2024 Fall   & 1 &  2 &  9 &  1 &  5 &  9 &  4 &  1 &  4 &  4 &  5 & 44 \\
                    & 2 &  3 & 15 &  0 &  5 &  6 &  6 &  0 &  5 &  1 &  1 & 42 \\
                    & 3 &  6 & 11 &  5 &  5 &  5 &  2 &  1 &  2 &  2 &  4 & 43 \\
                    & 4 &  3 &  4 &  3 &  6 &  4 &  6 &  3 &  4 &  3 &  6 & 42 \\
                    & 5 &  1 &  4 &  9 &  2 &  4 &  3 &  1 &  5 &  2 &  7 & 38 \\
                    & 6 &  2 &  2 &  7 &  4 &  6 &  4 &  2 &  2 &  0 &  5 & 34 \\
                    & 7 &  1 &  9 &  1 &  4 &  3 &  1 &  2 &  2 &  2 &  4 & 29 \\
        \addlinespace
        2025 Spring & 1 &  3 &  7 &  4 &  5 &  3 &  1 &  1 &  2 &  1 &  5 & 32 \\
                    & 2 &  5 &  2 &  5 &  2 &  3 &  5 &  1 &  5 &  1 &  7 & 36 \\
                    & 3 &  3 &  6 &  3 &  4 &  5 &  1 &  3 &  2 &  2 &  5 & 34 \\
                    & 4 &  5 &  3 &  6 &  1 &  1 &  1 &  3 &  3 &  1 &  8 & 32 \\
                    & 5 &  2 &  2 &  4 &  2 &  4 &  2 &  1 &  1 &  0 &  7 & 25 \\
        \addlinespace
        2025 Fall   & 1 &  7 & 18 &  8 &  3 &  1 &  2 &  0 &  0 &  2 &  3 & 44 \\
                    & 2 &  3 & 13 & 10 &  6 &  5 &  3 &  0 &  0 &  2 &  3 & 45 \\
                    & 3 &  3 &  8 &  9 &  5 &  5 &  4 &  1 &  2 &  0 &  5 & 42 \\
                    & 4 &  0 &  7 &  8 &  2 & 10 &  4 &  3 &  2 &  4 &  5 & 45 \\
        \bottomrule
    \end{tabular}
    \caption{Complete grade distributions by term and course section gathered from the UVA IRA (ira.virginia.edu/university-data-home). \textbf{N} denotes section enrollment. DFW combines D, F, and withdrawal outcomes. Treatment sections in 2024 Fall are 5, 6; in 2025 Spring just Sec.\ 1, and in 2025 Fall they are 1, 2.}
    \label{tab:grade_distributions_full}
\end{table}

\end{document}